# Variance and Covariance Inequalities for Truncated Joint Normal Distribution via Monotone Likelihood Ratio and Log-concavity

| Rahul Mukerjee | and | S.H. Ong |
|---|---|---|
| Indian Institute of Management Calcutta | | Institute of Mathematical Sciences |
| Joka, Diamond Harbour Road | | University of Malaya |
| Kolkata 700 104, India | | 50603 Kuala Lumpur, Malaysia |

*Abstract*: Let $X \sim N_v(0,\Lambda)$ be a normal vector in $v(\geq 1)$ dimensions, where $\Lambda$ is diagonal. With reference to the truncated distribution of $X$ on the interior of a $v$-dimensional Euclidean ball, we completely prove a variance inequality and a covariance inequality that were recently discussed by F. Palombi and S. Toti [*J. Multivariate Anal.* **122** (2013) 355-376]. These inequalities ensure the convergence of an algorithm for the reconstruction of $\Lambda$ only on the basis of the covariance matrix of $X$ truncated to the Euclidean ball. The concept of monotone likelihood ratio is useful in our proofs. Moreover, we also prove and utilize the fact that the cumulative distribution function of any positive linear combination of independent chi-square variates is log-concave, even though the same may not be true for the corresponding density function.

*Key words*: Chi-square distribution; covariance matrix reconstruction; positive linear combination; stochastic ordering.

*AMS Subject Classifications*: 60E15, 62E15.

## 1. Introduction

Let $X = (X_1,..., X_v)^T \sim N_v(0,\Lambda)$ be a normal vector in $v(\geq 1)$ dimensions, with zero means and covariance matrix $\Lambda = \mathrm{diag}(\lambda_1,...,\lambda_v)$. When truncated to the $v$-dimensional Euclidean ball $B(\rho)$ = $\{x : x^T x < \rho\}$, $X$ has density as given by

$$f(x) = k \prod_{n=1}^{v} \{\lambda_n^{-1/2} \phi(\lambda_n^{-1/2} x_n)\}, \text{ if } x = (x_1,..., x_v)^T \in B(\rho), \qquad (1.1)$$
$$= 0, \qquad \text{otherwise},$$

where $\phi(.)$ is the standard univariate normal density and $k$ is a normalizing constant. With reference to the truncated density (1.1), inspired by [4] and [3], we consider the inequalities

$$\mathrm{var}(X_n^2) \leq 2\lambda_n E(X_n^2),\ 1 \leq n \leq v, \qquad \mathrm{cov}(X_n^2, X_m^2) \leq 0, \quad 1 \leq n \neq m \leq v, \qquad (1.2)$$

which will hereafter be referred to as the variance inequality and the covariance inequality, respectively. While these authors stated the inequalities in (1.2) with respect to the conditional distribution of $X$ given $X \in B(\rho)$, our version in terms of the truncated distribution is equivalent to theirs and will somewhat simplify the notation in the sequel.

As noted in [3], the main interest in (1.2) originates from the fact that these inequalities, if universally true, are necessary and sufficient for the convergence of a fixed-point algorithm in [4] for the reconstruction of $\Lambda$ in case the only available information amounts to the covariance matrix of $X$ arising from the truncated density (1.1). Such reconstruction can be of practical importance, for example, in compositional analysis of multivariate log-normal data affected by outlying contaminations. Another motivation for (1.2) arises from non-linear optimization issues. We refer to [3] for further details. A complete proof of (1.2) is, however, quite nontrivial due to the symmetry mismatch between $N_v(0,\Lambda)$ and $B(\rho)$ which hinders exact calculation of moments. A very informative and in-depth discussion of the variance inequality was recently given in [3], where its truth was established for $0<\rho<2\lambda_n$ and $\rho >> \lambda_n$, and the intermediate region was left as an open challenge. The covariance inequality was not discussed in [3] and this also remains open.

The present paper aims at completely proving the two inequalities in (1.2). Following [3], we observe that, for any fixed *n*, the variance inequality is equivalent to

$$\partial E(\lambda_n^{-1} X_n^2)/\partial\lambda_n \le 0, \tag{1.3}$$

irrespective of $\lambda_m$, $m \ne n$. The concept of monotone likelihood ratio is useful in our proofs. Moreover, we also prove and utilize the fact that the cumulative distribution function (cdf) of any positive linear combination of independent chi-square variates is log-concave, even though the same may not be true for the corresponding density function. This log-concavity result, which is new to the best of our knowledge, should be of independent interest.

**2. Proof of variance inequality**

First let $v = 1$, which corresponds to the univariate case. Let $\Phi(.)$ be the standard univariate normal cdf and write $t=(\rho/\lambda_1)^{1/2}$. Then from (1.1), with the normalizing constant *k* there considered explicitly, it is not hard to see that $E(\lambda_1^{-1}X_1^2)=\mu(t)$, where $\mu(t)=1-2t\phi(t)/\{2\Phi(t)-1\}$. Note that

$$\mu'(t)=2\mu_1(t)\phi(t)/\{2\Phi(t)-1\}^2, \tag{2.1}$$

where the prime stands for differentiation and $\mu_1(t)=(t^2-1)\{2\Phi(t)-1\}+2t\phi(t)$. Since $\mu_1(t)$ tends to zero as $t\to 0+$ and $\mu_1'(t)=2t\{2\Phi(t)-1\}>0$ for $t > 0$, it follows that $\mu_1(t)\ge 0$ for $t > 0$, i.e., by (2.1), $\mu(t)$ is nondecreasing in *t* for $t > 0$. This implies that $E(\lambda_1^{-1}X_1^2)$ is nonincreasing in $\lambda_1$. Thus (1.3) and hence the variance inequality follow for $v = 1$.

Turing to the general case $v \ge 2$, we now prove the variance inequality for $n = 1$, without loss of generality. Let $Y_n=\lambda_n^{-1}X_n^2$, $1\le n\le v$. In the absence of truncation, $Y_1,\ldots,Y_v$ are independent, each distributed as chi-square with 1 degree of freedom (df). Hence with truncation as in (1.1), their joint density equals $k\Pi_{n=1}^{v} g_1(y_v)$, if $\Sigma_{n=1}^{v}\lambda_n y_n<\rho$, and 0 otherwise, where $g_j(.)$ is the chi-square density with *j* df. As a result, the marginal density of $Y_1$ is given by, say,

$$\psi(y_1;\lambda_1) = kg_1(y_1)H(\rho-\lambda_1 y_1), \quad \text{if } 0<y_1<\rho/\lambda_1, \tag{2.2}$$
$$= 0, \quad \text{otherwise,}$$

where, for $u > 0$, $H(u)$ is the integral of $\Pi_{n=2}^{v} g_1(y_v)$ over positive $y_2,\ldots,y_v$ satisfying $\Sigma_{n=2}^{v}\lambda_n y_n < u$. Clearly, $H(u)$ can be interpreted as the cdf of a linear combination of independent chi-square variates with positive coefficients. The $\lambda_1$ on the left-hand side of (2.2) makes the dependence of the marginal density on $\lambda_1$ explicit and is useful below in applying arguments based on MLR.

**Lemma 1**. *For any positive* $\lambda_{10},\lambda_{11}$, *satisfying* $\lambda_{10}<\lambda_{11}$, *the ratio* $\psi(y_1;\lambda_{11})/\psi(y_1;\lambda_{10})$ *is nonincreasing in* $y_1$ *over* $0<y_1<\rho/\lambda_{10}$.

*Proof.* By (2.2),

$$\psi(y_1;\lambda_{11})/\psi(y_1;\lambda_{10}) = k_0\{H(\rho-\lambda_{11}y_1)/H(\rho-\lambda_{10}y_1)\}, \quad \text{if } 0<y_1<\rho/\lambda_{11}, \tag{2.3}$$
$$= 0, \quad \text{if } \rho/\lambda_{11}\le y_1<\rho/\lambda_{10},$$

where $k_0$ is a positive constant which does not involve $y_1$. Differentiation with respect to $y_1$ shows that the ratio $H(\rho-\lambda_{11}y_1)/H(\rho-\lambda_{10}y_1)$ is nonincreasing in $y_1$ over $0<y_1<\rho/\lambda_{11}$ if and only if

$$\lambda_{10}\{H'(\rho-\lambda_{10}y_1)/H(\rho-\lambda_{10}y_1)\} \le \lambda_{11}\{H'(\rho-\lambda_{11}y_1)/H(\rho-\lambda_{11}y_1)\}. \tag{2.4}$$

Since $\lambda_{10}<\lambda_{11}$, (2.4) holds if $H'(u)/H(u)$ is nonincreasing in $u$, i.e., if $H(u)$ is log-concave for $u>0$. This, in turn, follows from Theorem 1 in Section 4 because, as mentioned above, $H(u)$ can be interpreted as the cdf of a linear combination of independent chi-square variates with positive coefficients. Therefore, (2.4) holds and the truth of the lemma is evident from (2.3). □

Since, with $\lambda_{10}<\lambda_{11}$, the support of $\psi(y_1;\lambda_1)$ is contained in $[0,\rho/\lambda_{10}]$ for both $\lambda_1=\lambda_{10}$ and $\lambda_1=\lambda_{11}$, the MLR property established in Lemma 1 induces a stochastic ordering of the marginal distributions of $Y_1$ for $\lambda_1=\lambda_{10}$ and $\lambda_1=\lambda_{11}$, and shows that the expectation of $Y_1$ for $\lambda_1=\lambda_{11}$ cannot exceed that for $\lambda_1=\lambda_{10}$, whenever $\lambda_{10}<\lambda_{11}$. Thus $E(Y_1)$ $[=E(\lambda_1^{-1}X_1^2)]$ is nonincreasing in $\lambda_1$, which proves (1.3) and hence the variance inequality for $v\ge 2$.

## 3. Proof of covariance inequality

Without loss of generality, we show that $\text{cov}(Y_1,Y_2)\le 0$, where $Y_n=\lambda_n^{-1}X_n^2$, $1\le n\le v$. For $u > 0$, let $\widetilde{H}(u)$ be the integral of $\Pi_{n=3}^{v} g_1(y_v)$ over positive $y_3,\ldots,y_v$ satisfying $\Sigma_{n=3}^{v}\lambda_n y_n < u$. Analogously to (2.2), then the marginal density of $(Y_1,Y_2)$ equals $kg_1(y_1)g_1(y_2)\widetilde{H}(\rho-\lambda_1 y_1-\lambda_2 y_2)$ if $y_1,y_2>0$ and $\lambda_1 y_1+\lambda_2 y_2<\rho$, and 0 otherwise. Hence by (2.2), the conditional density of $Y_2$ given $Y_1=y_1$, $0<y_1<\rho/\lambda_1$, turns out to be

$$\psi_{\text{cond}}(y_2 \mid y_1) = g_1(y_2)\widetilde{H}(\rho-\lambda_1 y_1-\lambda_2 y_2)/H(\rho-\lambda_1 y_1), \quad \text{if } 0<y_2<(\rho-\lambda_1 y_1)/\lambda_2, \tag{3.1}$$

$= 0,$ otherwise.

Note that as with $H(u)$ in Section 2, $\tilde{H}(u)$ is log-concave for $u > 0$, by Theorem 1 in Section 4.

**Lemma 2**. *For any fixed* $y_{10}, y_{11}$, *satisfying* $0 < y_{10} < y_{11} < \rho / \lambda_1$, *the ratio*

$$\psi_{\text{cond}}(y_2 \mid y_{11}) / \psi_{\text{cond}}(y_2 \mid y_{10})$$

*is nonincreasing in* $y_2$ *over* $0 < y_2 < (\rho - \lambda_1 y_{10}) / \lambda_2$.

*Proof.* By (3.1), the ratio considered here equals $\tilde{k}\{\tilde{H}(\rho - \lambda_1 y_{11} - \lambda_2 y_2) / \tilde{H}(\rho - \lambda_1 y_{10} - \lambda_2 y_2)\}$ if $0 < y_2 < (\rho - \lambda_1 y_{11}) / \lambda_2$, and 0 if $(\rho - \lambda_1 y_{11}) / \lambda_2 \le y_2 < (\rho - \lambda_1 y_{10}) / \lambda_2$. Here $\tilde{k}$ is a positive quantity which does not involve $y_2$. For $v = 2$, this ratio simply equals $\tilde{k}$ over $0 < y_2 < (\rho - \lambda_1 y_{11}) / \lambda_2$, because then the term involving $\tilde{H}(.)$ does not arise in (3.1); as a result, the lemma is immediate. On the other hand, for $v \ge 3$, this ratio is nonincreasing in $y_2$ over $0 < y_2 < (\rho - \lambda_1 y_{11}) / \lambda_2$, a fact which can be established as in Lemma 1 if one notes that $y_{10} < y_{11}$ and invokes the log-concavity of $\tilde{H}(u)$ for $u > 0$. Thus the lemma holds again. □

Analogously to Section 2, the MLR property shown in Lemma 2 implies that $E(Y_2 \mid Y_1 = y_1)$ is nonincreasing in $y_1$ over $0 < y_1 < \rho / \lambda_1$. Therefore, $\text{cov}\{Y_1, E(Y_2 \mid Y_1)\} \le 0$, and hence a conditioning argument yields $\text{cov}(Y_1, Y_2) \le 0$, proving the covariance inequality.

## 4. A log-concavity theorem

**Theorem 1.** *The cdf of any linear combination of independent chi-square variates with positive coefficients is log-concave over* $(0, \infty)$.

*Proof.* It will be convenient to present the proof through several steps.

Step 1 (Preliminaries): Due to the reproductive property of independent chi-squares, it suffices to prove the theorem for the case where each chi-square variate in the linear combination has 1 df. To that effect, let $U = \Sigma_{i=1}^{s} a_i Y^{(i)}$, where $a_1, ..., a_s$ are positive constants and $Y^{(1)}, ..., Y^{(s)}$ are independent chi-square variates each with 1 df. Let $H(.)$ be the cdf of $U$. Without loss of generality, suppose

$$a_1 \le ... \le a_s \quad \text{and} \quad a_1 = 1. \tag{4.1}$$

While the first condition in (4.1) is clearly allowable, the second condition also entails no loss of generality as the log-concavity of $H(u)$ over $u > 0$ is equivalent to that of $H(u / a)$, for any $a > 0$.

Step 2 (A background result): If $a_1 = ... = a_s = 1$, then the theorem follows from the well-known log-concavity of the cdf of a chi-square variate; see e.g., [1]. Let, therefore, $s \ge 2$ and suppose at least one of $a_2, ..., a_s$ is greater than 1. In view of (4.1), then following Theorem 2 in [6],

$$H(u) = \sum_{j=0}^{\infty} p_j G_{s+2j}(u), \tag{4.2}$$

where $G_{s+2j}(.)$ is the chi-square cdf with $s+2j$ df,

$$p_0 = \prod_{i=2}^{s} a_i^{-1/2} \qquad \text{and} \qquad p_j = j^{-1} \sum_{i=0}^{j-1} M_{j-i} p_i , \quad j \geq 1, \tag{4.3}$$

with

$$M_j = \tfrac{1}{2}(c_2^j + ... + c_s^j), \quad j \geq 1, \quad \text{and} \quad c_i = 1 - a_i^{-1}, \quad 2 \leq i \leq s . \tag{4.4}$$

Incidentally, (4.2) was reported earlier in [5] without, however, the recursion formula (4.3) that will be crucial to us. From (4.1) and (4.4), $0 \leq c_i < 1$ for each $i$, and at least one of them is positive, since at least one of $a_2,...,a_s$ exceeds 1. By (4.3), this implies that each $p_j$ is positive. Indeed, as noted in [6, p. 545], $\Sigma_{j=0}^{\infty} p_j = 1$, i.e., $\{p_j\}$ is a probability sequence. Consequently, the sequence of partial sums $\{W_i\}$, where $W_i = \Sigma_{j=0}^{i} p_j$, $i \geq 0$, is bounded.

Step 3 (Power series expansion): By standard results on power series expansion for the chi-square cdf, from (4.2), we now get, for $u > 0$,

$$\begin{aligned} H(u) &= \exp(-u/2) \sum_{j=0}^{\infty} \sum_{i=j}^{\infty} p_j (u/2)^{s/2+i} / \Gamma(s/2+i+1) \\ &= \exp(-u/2) \sum_{i=0}^{\infty} W_i (u/2)^{s/2+i} / \Gamma(s/2+i+1) \\ &= \exp(-u/2)(u/2)^{s/2-1} \beta(u) \end{aligned} \tag{4.5}$$

where

$$\beta(u) = \sum_{i=0}^{\infty} W_i (u/2)^{i+1} / \Gamma(s/2+i+1) = \sum_{i=0}^{\infty} W_{i-1} (u/2)^i / \Gamma(s/2+i), \tag{4.6}$$

with $W_{-1} = 0$. The change in the order of summation to reach (4.5) is valid because $\{W_i\}$ is bounded and hence the power series $\beta(u)$ has infinite radius of convergence. The point just noted also validates term-by-term differentiation of $\beta(u)$ as used in Steps 4 and 6 below.

Step 4 (Proof for $s = 2$): Let $s = 2$. Then from [5, p. 556],

$$p_j = p_0 \binom{2j}{j} (c_2/4)^j, \quad j \geq 0. \tag{4.7}$$

As $c_2 < 1$, this yields $p_{j+1}/p_j = (2j+1)c_2/(2j+2) < 1$, $j \geq 0$, i.e., $\{p_j\}$ is a decreasing sequence. From (4.5) and (4.6) on simplification, it now follows that for $u > 0$,

$$\begin{aligned} H''(u) &= \tfrac{1}{4} \exp(-u/2) \sum_{i=0}^{\infty} (W_{i+1} - 2W_i + W_{i-1})(u/2)^i / i! \\ &= \tfrac{1}{4} \exp(-u/2) \sum_{i=0}^{\infty} (p_{i+1} - p_i)(u/2)^i / i! \ \leq 0 . \end{aligned}$$

Thus $H(u)$ is concave, and hence log-concave, over $(0, \infty)$ and the theorem is proved for $s = 2$.

Step 5 (Discussion on $s \geq 3$) The approach for $s = 2$ does not work for general $s$, due to lack of the nonincreasing property of $\{p_j\}$; e.g., if $s = 4$ and $a_1 = 1$, $a_2 = a_3 = a_4 = 5$, then by (4.4), $c_2 = c_3 = c_4 = 0.8$, $M_1 = 1.2$, and (4.3) yields $p_1 = 1.2 p_0$. For $s \geq 3$, however, $s/2 > 1$ and in view of (4.5) and (4.6), one may wonder if the log-concavity of $H(u)$ can be deduced from Theorem 2.1 in [2, p. 107]. This would require (a) $\Sigma_{i=0}^{\infty} W_i < \infty$, and (b) log-concavity of the sequence $\{W_i\}$. While Lemma 3 below proves (b), condition (a) is not met here because $W_i$ tends to $\Sigma_{j=0}^{\infty} p_j$ (=1) as $i \to \infty$. Nevertheless, the arguments in [2] go through, because $\beta(u)$ has infinite radius of convergence and hence allows term-by-term differentiability. In order to convince the reader and also for completeness, the details are shown in the next step.

Step 6 (Proof for $s \geq 3$) Let $s \geq 3$. Since $\exp(-u/2)$ is log-concave, in view of (4.5) it suffices to show the log-concavity of $L(u) = (u/2)^{s/2-1}\beta(u)$ over $(0,\infty)$. As $W_{-1} = 0$, by (4.6), for $u > 0$.

$$L'(u) = \tfrac{1}{2}(u/2)^{s/2-2} \sum_{i=0}^{\infty} W_i (u/2)^{i+1} / \Gamma(s/2+i) = \tfrac{1}{2}(u/2)^{s/2-2} \sum_{i=0}^{\infty} W_{i-1} (u/2)^{i} / \Gamma(s/2+i-1),$$

$$L''(u) = \tfrac{1}{4}(u/2)^{s/2-2} \sum_{i=0}^{\infty} W_i (u/2)^{i} / \Gamma(s/2+i-1).$$

Note that the case $s = 2$ required a separate treatment because then $\Gamma(s/2+i-1)$ is undefined at $i = 0$ and hence the aforesaid expression for $L''(u)$ is invalid. For $i \geq 0$, write

$$d(0,i) = W_i, \qquad d(1,i) = W_{i-1}, \tag{4.8}$$

$$q(i,0) = (u/2)^i / \Gamma(s/2+i-1), \qquad q(i,1) = (u/2)^{i+1} / \Gamma(s/2+i). \tag{4.9}$$

Also, let $Z = (z_{jm})$ and $Z^{(n)} = (z_{jm}^{(n)})$, $n \geq 1$, be $2 \times 2$ matrices such that $z_{jm} = \Sigma_{i=0}^{\infty} d(j,i) q(i,m)$ and $z_{jm}^{(n)} = \Sigma_{i=0}^{n} d(j,i) q(i,m)$, for $j, m = 0,1$. From (4.6), (4.8), (4.9), and the expressions for $L'(u)$ and $L''(u)$ as shown above, then for $u > 0$,

$$L''(u) = \tfrac{1}{4}(u/2)^{s/2-2} z_{00}, \quad L'(u) = \tfrac{1}{2}(u/2)^{s/2-2} z_{01} = \tfrac{1}{2}(u/2)^{s/2-2} z_{10}, \quad L(u) = (u/2)^{s/2-2} z_{11},$$

so that

$$L(u)L''(u) - \{L'(u)\}^2 = \tfrac{1}{4}(u/2)^{s-4} \det(Z) = \tfrac{1}{4}(u/2)^{s-4} \lim_{n\to\infty} \det(Z^{(n)}). \tag{4.10}$$

Now, $Z^{(n)} = D_n Q_n$, where $D_n$ is a $2 \times (n+1)$ matrix with $(j,i)$th element $d(j,i)$ $(j = 0,1; 0 \leq i \leq n)$, and $Q_n$ is an $(n+1) \times 2$ matrix with $(i,m)$th element $q(i,m)$ $(0 \leq i \leq n; m = 0,1)$. For $0 \leq i < r \leq n$, let $D_n(i,r)$ be the $2 \times 2$ submatrix of $D_n$ consisting of its $i$th and $r$th columns, and $Q_n(i,r)$ be the $2 \times 2$ submatrix of $Q_n$ consisting of its $i$th and $r$th rows. Then by the Cauchy-Binet formula,

$$\det(Z^{(n)}) = \det(D_n Q_n) = \sum_{i=0}^{n-1} \sum_{r=i+1}^{n} \det\{D_n(i,r)\} \det\{Q_n(i,r)\}. \tag{4.11}$$

For $0 \leq i < r \leq n$, by (4.8), Lemma 3 below and the fact that $W_{-1} = 0$, we have $\det\{D_n(i,r)\}$ = $W_i W_{r-1} - W_r W_{i-1} \geq 0$, while after some simplification, (4.9) yields

$$\det\{Q_n(i,r)\} = (u/2)^{i+r+1}(i-r) / \{\Gamma(s/2+i)\Gamma(s/2+r)\} < 0.$$

Hence by (4.11), $\det(Z^{(n)}) \le 0$ for every $n \ge 1$. Therefore, $\lim_{n\to\infty} \det(Z^{(n)}) \le 0$ and the log-concavity of $L(u)$ is immediate from (4.10). This completes the proof of the theorem. □

**Remark 1**. In view of the similarity between (4.2) and a noncentral chi-square cdf, one may wonder if Theorem 1 could be proved along the lines of [1]. This possibility is precluded by the fact that, unlike the Poisson weights appearing in a noncentral chi-square cdf, the sequence $\{p_j\}$ here is not log-concave in general; e.g., with $s = 2$, by (4.7), $p_1 = \frac{1}{2}c_2 p_0$, $p_2 = \frac{3}{8}c_2^2 p_0$, and $p_1^2 < p_0 p_2$, as $c_2 > 0$. Indeed, unlike with a noncentral chi-square, the density $h(u)$ corresponding to our $H(u)$ may not be log-concave. Thus, with $s = 2$ again, from (4.5) one can check that $h(u)h''(u) - \{h'(u)\}^2$ tends to $(p_0 p_2 - p_1^2)/16 > 0$ as $u \to 0+$. □

**Lemma 3**. *The sequence* $\{W_i\}$ *is log-concave, i.e.,* $W_i^2 \ge W_{i-1}W_{i+1}$, $i \ge 1$.

*Proof.* Since $W_i = \Sigma_{j=0}^{i} p_j$, one can check that $W_i^2 \ge W_{i-1}W_{i+1}$ if and only if $p_{i+1}/p_i \le W_i/W_{i-1}$. It, therefore, suffices to prove the inequality

$$\{(i+1)p_{i+1}\}/(ip_i) \le W_i/W_{i-1}, \quad i \ge 1, \tag{4.12}$$

which is even stronger. Since $0 \le c_i < 1$, for $2 \le i \le s$, it is clear from (4.4) that

$$M_j \ge M_{j+1}, \quad j \ge 1. \tag{4.13}$$

Also, from (4.3), $p_1 = M_1 p_0$ and $2p_2 = M_2 p_0 + M_1 p_1 = (M_2 + M_1^2)p_0$. Hence by (4.13), $2p_2/p_1$ = $M_2 M_1^{-1} + M_1 \le 1 + M_1 = (p_0 + p_1)/p_0 = W_1/W_0$, i.e., (4.12) holds for $i = 1$. To apply the method of induction, let (4.12) hold for $1 \le i \le n$. If possible, suppose (4.12) does not hold for $i = n+1$. Then by (4.3), $\Sigma_{i=0}^{n+1} M_{n+2-i} p_i / \Sigma_{i=0}^{n} M_{n+1-i} p_i > W_{n+1}/W_n$, i.e., recalling (4.4),

$$\frac{\sum\limits_{j=2}^{s}\sum\limits_{i=0}^{n+1} c_j^{n+2-i} p_i}{\sum\limits_{j=2}^{s}\sum\limits_{i=0}^{n} c_j^{n+1-i} p_i} > \frac{W_{n+1}}{W_n},$$

which implies that $\Sigma_{i=0}^{n+1} c_j^{n+2-i} p_i / \Sigma_{i=0}^{n} c_j^{n+1-i} p_i > W_{n+1}/W_n$, for some *j*. Obviously, $c_j > 0$ for this *j*. Moreover, as $c_j < 1$, writing *c* for this $c_j$, it follows that there exists *c*, $0 < c < 1$, such that

$$\frac{\sum\limits_{i=0}^{n+1} c^{n+2-i} p_i}{\sum\limits_{i=0}^{n} c^{n+1-i} p_i} > \frac{W_{n+1}}{W_n} = \frac{\sum\limits_{i=0}^{n+1} p_i}{\sum\limits_{i=0}^{n} p_i}, \quad \text{i.e.,} \quad c + \frac{p_{n+1}}{\sum\limits_{i=0}^{n} c^{n-i} p_i} > 1 + \frac{p_{n+1}}{\sum\limits_{i=0}^{n} p_i},$$

$$\text{i.e.,} \qquad (1-c)(\sum_{i=0}^{n} c^{n-i} p_i)(\sum_{i=0}^{n} p_i) < p_{n+1} \sum_{i=0}^{n-1} (1-c^{n-i}) p_i,$$

i.e., $$(\sum_{i=0}^{n} c^{n-i} p_i)(\sum_{i=0}^{n} p_i) < p_{n+1} \sum_{i=0}^{n-1} \sum_{j=0}^{n-i-1} c^j p_i \,, \tag{4.14}$$

division of both sides by $1-c$ in the last step being permissible as $0<c<1$.

For $n=1$, (4.14) reduces to $(cp_0+p_1)(p_0+p_1)<p_2p_0$, and as $cp_0 > 0$, this yields $p_2p_0 > p_1(p_0+p_1)$, i.e., $p_2/p_1 > W_1/W_0$. As a result, $2p_2/p_1 > W_1/W_0$, which violates the truth of (4.12) for $i = 1$.

Next, suppose $n \ge 2$. Since by induction hypothesis, (4.12) holds for $i=n$, we get $p_{n+1}/p_n < \{(n+1)p_{n+1}\}/(np_n) \le W_n/W_{n-1}$, i.e., $p_{n+1} < p_n \Sigma_{i=0}^{n} p_i / \Sigma_{i=0}^{n-1} p_i$. Using this in (4.14),

$$(\sum_{i=0}^{n} c^{n-i} p_i)(\sum_{i=0}^{n-1} p_i) < p_n \sum_{i=0}^{n-1} \sum_{j=0}^{n-i-1} c^j p_i \,,$$

i.e., $$(\sum_{i=0}^{n-1} c^{n-i} p_i + p_n)(\sum_{i=0}^{n-1} p_i) < p_n(\sum_{i=0}^{n-2} \sum_{j=1}^{n-i-1} c^j p_i + \sum_{i=0}^{n-2} p_i + p_{n-1}) \,,$$

i.e., $$(\sum_{i=0}^{n-1} c^{n-i} p_i)(\sum_{i=0}^{n-1} p_i) < p_n(\sum_{i=0}^{n-2} \sum_{j=1}^{n-i-1} c^j p_i),$$

i.e., $$(\sum_{i=0}^{n-1} c^{n-i-1} p_i)(\sum_{i=0}^{n-1} p_i) < p_n(\sum_{i=0}^{n-2} \sum_{j=1}^{n-i-1} c^{j-1} p_i) = p_n(\sum_{i=0}^{n-2} \sum_{j=0}^{n-i-2} c^j p_i), \tag{4.15}$$

as $c > 0$. Note that (4.15) has the same form as (4.14) with $n$ in (4.14) replaced by $n-1$. Continuation of the above steps, with the use of the induction hypothesis for $i = n, n–1,\ldots,$ eventually leads to (4.14) with $n = 1$. But, as noted in the last paragraph, this contradicts the truth of (4.12) for $i = 1$.

Thus if (4.12) holds for $1 \le i \le n$, then it must hold for $i = n+1$. Since (4.12) holds for $i = 1$, this establishes the truth of (4.12) for every $i \ge 1$ and hence proves the lemma. □

**Acknowledgement**. The work of Mukerjee was supported by the J.C. Bose National Fellowship of the Government of India and a grant from the Indian Institute of Management Calcutta. The work of Ong was supported by the University of Malaya UMRGS grant RP009A-13AFR.